\documentclass[12pt,a4paper]{amsart}
\usepackage{mathrsfs}
\usepackage{xcolor}
\usepackage{graphicx}
\newtheorem{theo+}              {Theorem}           [section]
\newtheorem{prop+}  [theo+]     {Proposition}
\newtheorem{coro+}  [theo+]     {Corollary}
\newtheorem{lemm+}  [theo+]     {Lemma}
\newtheorem{exam+}  [theo+]     {Example}
\newtheorem{rema+}  [theo+]     {Remark}
\newtheorem{defi+}  [theo+]     {Definition}

\newenvironment{theorem}{\begin{theo+}}{\end{theo+}}
\newenvironment{proposition}{\begin{prop+}}{\end{prop+}}
\newenvironment{corollary}{\begin{coro+}}{\end{coro+}}

\newenvironment{definition}{\begin{defi+}}{\end{defi+}}

\usepackage{amsthm}
\theoremstyle{plain} \theoremstyle{remark}
\newtheorem{remark}{Remark}
\newtheorem{example}{Example}

\def \r{\mbox{${\mathbb R}$}}
\def\E{/\kern-1.0em \equiv }

\title[On the stress-bienergy tensor of conformal immersions]
{On the stress-bienergy tensor of conformal immersions and hypersurfaces}

\author{A. Mohammed Cherif}
\address{Faculty of Exact Sciences,\newline\indent University Mustapha Stambouli,
\newline\indent Mascara, Algeria.\newline\indent Email: a.mohammedcherif@univ-mascara.dz}

\begin{document}
\maketitle

\begin{abstract}
We study the stress-bienergy tensor $\bar S_2$ of a conformal immersion
$\phi:(M^m,\bar g)\to(N^n,h)$ with conformal factor $\lambda$, and we compute it,
together with its trace and its divergence, in terms of $\lambda$ and of the
second fundamental form of the associated isometric immersion. This gives a
characterization of biconservative conformal immersions and, for hypersurfaces,
we prove that $\mathrm{div}\,\bar S_2=0$ if and only if the tangent part of the
bitension field vanishes. We introduce conformally biconservative hypersurfaces
and show that a totally umbilical hypersurface of a space form is conformally
biconservative with a non-constant conformal factor if and only if this factor
is an isoparametric function satisfying an explicit ordinary differential
equation. We illustrate this result with several examples, in which dimension
four plays a special role. Finally, we study the equation $\bar S_2=0$: on a
compact manifold of dimension $m\neq 4$, it forces the associated isometric
immersion to be harmonic.
\end{abstract}

\section{Introduction}

Let $\phi:(M^m,g)\to(N^n,h)$ be a smooth map between Riemannian manifolds.
Harmonic maps \cite{ES,EL1,EL2} are the critical points of the energy
\[
E(\phi)=\frac12\int_M|d\phi|^2\,dv_g,
\]
with Euler--Lagrange equation $\tau(\phi)=\mathrm{Tr}_g\nabla d\phi=0$. Baird and
Eells \cite{BE} associated to $E$ the stress-energy tensor
$S_\phi=e(\phi)\,g-\phi^*h$ (see also \cite{B2}), where $e(\phi)=\frac12|d\phi|^2$,
which satisfies $\mathrm{div}\,S_\phi=-h(\tau(\phi),d\phi)$. Biharmonic maps,
introduced by Jiang \cite{J}, are the critical points of the bienergy
\[
E_2(\phi)=\frac12\int_M|\tau(\phi)|^2\,dv_g;
\]
see \cite{O1,Ou6,OC} for background and surveys. The analogue of $S_\phi$ for
$E_2$ is the stress-bienergy tensor $S_2$, studied by Loubeau, Montaldo and
Oniciuc \cite{LMO}. It satisfies
$\mathrm{div}\,S_2=-\langle\tau_2(\phi),d\phi\rangle$, where $\tau_2(\phi)$ is
the bitension field. A map is called biconservative if
$\mathrm{div}\,S_2=0$. Every biharmonic map is biconservative, but the converse
fails in general. For an isometric immersion, biconservativity means that the
tangent part of $\tau_2(\phi)$ vanishes.\\
In this paper we study $S_2$ for conformal immersions. Let
$\phi:(M^m,\bar g)\to(N^n,h)$ be an immersion with $\phi^*h=\lambda^2\bar g$
for a smooth function $\lambda:M\to\mathbb{R}^+$. Then $\phi$ is also an
isometric immersion of $(M^m,g)$ into $(N^n,h)$, with $g=\lambda^2\bar g$.
Conformal biharmonic immersions were investigated in
\cite{BFO,LO2,Ou1,Ou2,Ou4,BO,Ou5,WC,MO}. Here we ask how the stress-bienergy
tensor $\bar S_2$ of $\phi:(M,\bar g)\to(N,h)$ depends on the geometry of the
associated isometric immersion (its second fundamental form $B$ and mean
curvature vector $\eta$) and on the conformal factor $\lambda$.\\
Our main results are as follows. We give explicit formulas for $\bar S_2$
(Theorem \ref{th1}), for its trace (Proposition \ref{prop1}) and for its
divergence (Theorem \ref{th2}). This yields a characterization of biconservative
conformal immersions. For hypersurfaces, we show that
$\mathrm{div}_{\bar g}\bar S_2=0$ if and only if the tangent part of the
bitension field $\tau_2(\phi,\bar g)$ vanishes (Theorem \ref{th3}); compare
\cite{OuP} for biharmonic hypersurfaces. We then introduce conformally
biconservative hypersurfaces. For totally umbilical hypersurfaces of space forms,
we prove that a non-constant conformal factor yields a conformally biconservative
hypersurface if and only if it is an isoparametric function whose structure
functions satisfy an explicit ODE (Theorem \ref{thm3}); this is the
biconservative counterpart of the isoparametric property of the conformal factor
for biharmonic conformal immersions obtained in \cite{MO}. We illustrate this
with several examples, in which the dimension $m=4$ plays a special role.
Finally, we consider the equation $\bar S_2=0$ and show that, when $m\neq 4$
and $M$ is compact, it forces the associated isometric immersion to be harmonic.
In the non-compact case, we construct non-trivial solutions.

%%%%%%%%%%%%%%%%%%%%%%%%%%%%%%%%%%%%%%%%%%%%%%%%%%%%%%%%%%%%%%%%%%%%%%%%%%%%%%%%%%%%%%%%%%%%%%%%%%%%%%%%%%%%%%%%%%%%%%%%%%%%%%%%%%%%%%%%%%%%

\section{The stress-bienergy tensor of conformal immersions}

In this section, we express $\bar S_2$, its trace and its divergence in terms of
the conformal factor $\lambda$ and of the second fundamental form $B$ of the
associated isometric immersion. We deduce a criterion for biconservativity and a
condition under which $\bar S_2=0$ forces harmonicity.

\begin{theorem}\label{th1}
Let $\phi:(M^m,\bar{g})\to (N^n, h)$ be a conformal immersion with conformal factor $\lambda: M^m\to \r^+$.
Then, the biharmonic stress tensor of the conformal immersion $\phi$ is given by
\begin{eqnarray*}
 \bar{S}_2(X,Y)
   &=& \frac{1}{2}\Big[(m-2)(m-4)|\nabla\lambda|^2-2(m-2)\lambda \Delta\lambda-m^2\lambda^2|\eta|^2\Big]g(X,Y)
       \\
   &&+2m\lambda^2h(B(X,Y),\eta)
     +2(m-2)\left[X(\lambda)Y(\lambda)+\lambda ({\rm Hess\,}\lambda)(X,Y)\right],
\end{eqnarray*}
for all $X,Y\in\Gamma(TM)$, where $\eta$ and $B$ are the mean curvature vector, the second fundamental form of the associated isometric immersion $\phi:(M^m,g)\to (N^n, h)$ respectively, where $g=\phi^{*}h=\lambda^2 {\bar g}$, and the operators $\nabla$, ${\rm Hess}$, $\Delta$  and $|,|$ are taken with respect to $g$.
\end{theorem}

\begin{proof}
We recall that
\begin{eqnarray}\label{th1-eq1}
 \bar{S}_2(X,Y)
   &=&\nonumber \frac{1}{2}|\tau(\phi,\bar{g})|^2\bar{g}(X,Y)+\langle {\rm d}\phi,\nabla \tau(\phi,\bar{g})\rangle \bar{g}(X,Y)\\
   & & -h( {\rm d}\phi(X),\nabla_Y^\phi \tau(\phi,\bar{g}))- h( {\rm d}\phi(Y),\nabla_X^\phi \tau(\phi,\bar{g})),
\end{eqnarray}
for all $X,Y\in\Gamma(TM)$ (see \cite{LMO}). Note that the tension field of the conformal immersion $\phi$ is given by
\begin{equation}\label{th1-eq2}
    \tau(\phi,\bar{g})=m\lambda^2\eta-(m-2)\lambda {\rm d}\phi(\nabla\lambda),
\end{equation}
(see \cite{Ou1}). Let $X,Y\in\Gamma(TM)$. Using $h(\eta,{\rm d}\phi(X))=0$ and $g=\phi^{*}h$, we get
\begin{eqnarray}\label{th1-eq3}
 |\tau(\phi,\bar{g})|^2 \bar{g}(X,Y)
   &=&\nonumber h(\tau(\phi,\bar{g}),\tau(\phi,\bar{g}))\lambda^{-2}g(X,Y) \\
   &=&\nonumber \left(m^2\lambda^4|\eta|^2+(m-2)^2\lambda^2|{\rm d}\phi(\nabla\lambda)|^2\right) \lambda^{-2}g(X,Y) \\
   &=& \left(m^2\lambda^2|\eta|^2+(m-2)^2|\nabla\lambda|^2\right)g(X,Y),
\end{eqnarray}
and the following
\begin{eqnarray}\label{th1-eq4}
h({\rm d}\phi(X),\nabla_{Y}^\phi \tau(\phi,\bar{g}))
   &=&\nonumber Y(h({\rm d}\phi(X), \tau(\phi,\bar{g})))-h(\nabla_{Y}^\phi{\rm d}\phi(X), \tau(\phi,\bar{g})) \\
   &=&\nonumber -(m-2)Y(\lambda h({\rm d}\phi(X), {\rm d}\phi(\nabla\lambda)))-m\lambda^2h(\nabla_{Y}^\phi{\rm d}\phi(X),\eta ) \\
   & &\nonumber +(m-2)\lambda h(\nabla_{Y}^\phi{\rm d}\phi(X), {\rm d}\phi(\nabla\lambda) ) \\
   &=&\nonumber -(m-2)Y(\lambda g(X,\nabla\lambda))-m\lambda^2h(B(Y,X),\eta ) \\
   & &\nonumber +(m-2)\lambda h({\rm d}\phi(\nabla_{Y}^MX), {\rm d}\phi(\nabla\lambda) ) \\
   &=&\nonumber -(m-2)Y(\lambda) X(\lambda) -(m-2)\lambda g(\nabla_{Y}^MX,\nabla\lambda)\\
   & &\nonumber -(m-2)\lambda g(X,\nabla_{Y}^M\nabla\lambda)-m\lambda^2h(B(Y,X),\eta ) \\
   & &\nonumber +(m-2)\lambda g(\nabla_{Y}^MX,\nabla\lambda)\\
   &=&\nonumber -(m-2)Y(\lambda) X(\lambda)-(m-2)\lambda ({\rm Hess\,}\lambda)(Y,X) \\
   & &          -m\lambda^2h(B(Y,X),\eta ).
\end{eqnarray}
Let $\{\bar{e}_i\}_{1\leq i\leq m}$ be a local orthonormal frame in $M^m$ with respect to the Riemannian metric $\bar{g}$. Then, $\{e_i=\frac{1}{\lambda}\bar{e}_i\}_{1\leq i\leq m}$ is a local orthonormal frame in $M^m$ with respect to the Riemannian metric $g=\lambda^2\bar{g}$. Using (\ref{th1-eq4}), we obtain
\begin{eqnarray}\label{th1-eq5}
\langle d\phi,\nabla \tau(\phi,\bar{g})\rangle
   &=&\nonumber \sum_{i=1}^m h(d\phi(\bar{e}_i),\nabla^\phi_{\bar{e}_i} \tau(\phi,\bar{g})) \\
   &=&\nonumber \lambda^2\sum_{i=1}^m h(d\phi({e}_i),\nabla^\phi_{{e}_i} \tau(\phi,\bar{g})) \\
   &=&\nonumber \lambda^2\sum_{i=1}^m\Big[-(m-2)e_i(\lambda) e_i(\lambda)-(m-2)\lambda ({\rm Hess\,}\lambda)(e_i,e_i) \\
   & &          -m\lambda^2h(B(e_i,e_i),\eta )\Big].
\end{eqnarray}
Therefore, it follows from (\ref{th1-eq5}) that
\begin{eqnarray}\label{th1-eq6}
\langle d\phi,\nabla \tau(\phi,\bar{g})\rangle
   &=& -\lambda^2\left[(m-2)|\nabla\lambda|^2+(m-2)\lambda \Delta\lambda+m^2\lambda^2|\eta|^2\right].
\end{eqnarray}
The proof of Theorem \ref{th1} follows from equations \eqref{th1-eq1}, \eqref{th1-eq3}, \eqref{th1-eq4}, and \eqref{th1-eq6}.
\end{proof}

\begin{remark}\label{remark1}
i) For $\lambda=1$, the tensor $\bar S_2$ becomes
\begin{eqnarray*}
 S_2(X,Y)
   &=& -\frac{1}{2}m^2|\eta|^2g(X,Y)+2m h(B(X,Y),\eta).
\end{eqnarray*}
By using $\tau(\phi,g)=m\eta$, we find that (see \cite{LMO})
\begin{eqnarray}\label{re1-eq1}
 S_2(X,Y)
   &=& -\frac{1}{2}|\tau(\phi,g)|^2g(X,Y)+2h(\tau(\phi,g),B(X,Y)).
\end{eqnarray}
ii) For $m=2$, Theorem \ref{th1} gives $\bar S_2=\lambda^2 S_2$, where $S_2$ is the
stress-bienergy tensor of the isometric immersion $\phi:(M^2,g)\to(N^n,h)$.\\
iii) When $n=m$, a conformal map $\phi:(M^m,\bar g)\to(N^m,h)$ is a local diffeomorphism, so the associated isometric immersion $\phi:(M,g)\to(N,h)$ is a local isometry. Hence $B=0$ and $\eta=0$, and Theorem \ref{th1} reduces to
\begin{eqnarray*}
\bar S_2(X,Y)&=&\frac12(m-2)\Big[(m-4)|\nabla\lambda|^2-2\lambda\Delta\lambda\Big]g(X,Y)\\
&&+2(m-2)\Big[X(\lambda)Y(\lambda)+\lambda\,(\operatorname{Hess}\lambda)(X,Y)\Big].
\end{eqnarray*}
In particular, $\bar S_2$ depends only on the conformal factor $\lambda$, which agrees with the result of \cite{OD} for conformal maps between manifolds of the same dimension. Formula \eqref{th1-eq2} gives $\tau(\phi,\bar g)=-(m-2)\lambda\, d\phi(\nabla\lambda)$.
\end{remark}

\begin{proposition}\label{prop1}
Let $\phi:(M^m,\bar{g})\to (N^n, h)$ be a conformal immersion with conformal factor $\lambda: M^m\to \r^+$.
Then, the trace of the biharmonic stress tensor of $\phi$ is given by
\begin{eqnarray*}
{\rm Tr}_{\bar{g}}\,\bar{S}_2
&=&\frac{1}{2}\lambda^2 \Big[
   (m-2)^3|\nabla\lambda|^2
   -2(m-2)^2 \lambda \Delta\lambda
   - m^2 (m-4)\lambda^2|\eta|^2
   \Big].
\end{eqnarray*}
\end{proposition}

\begin{proof}
Let $\{\bar{e}_i\}_{1\leq i\leq m}$ be a local orthonormal frame in $(M^m,\bar{g})$. We have, $\{e_i=\frac{1}{\lambda}\bar{e}_i\}_{1\leq i\leq m}$ is a local orthonormal frame in $(M^m,g)$. According to Theorem \ref{th1}, we find that
\begin{eqnarray}\label{prop1-eq7}
{\rm Tr}_{\bar{g}}\,\bar{S}_2
&=&\nonumber \sum_{i=1}^m \bar{S}_2(\bar{e}_i,\bar{e}_i)\\
&=&\nonumber \lambda^2\sum_{i=1}^m \bar{S}_2({e}_i,{e}_i)\\
&=&\nonumber \lambda^2\sum_{i=1}^m \Big\{
\frac{1}{2}\Big[(m-2)(m-4)|\nabla\lambda|^2-2(m-2)\lambda \Delta\lambda-m^2\lambda^2|\eta|^2\Big]g(e_i,e_i)\\
&&+2m\lambda^2h(B(e_i,e_i),\eta)
+2(m-2)\left[e_i(\lambda)e_i(\lambda)+\lambda ({\rm Hess\,}\lambda)(e_i,e_i)\right]\Big\}.
\end{eqnarray}
Therefore, it follows from (\ref{prop1-eq7}) that
\begin{eqnarray}\label{prop1-eq8}
{\rm Tr}_{\bar{g}}\,\bar{S}_2
&=&\nonumber \lambda^2 \Big\{
\frac{m}{2}\Big[(m-2)(m-4)|\nabla\lambda|^2-2(m-2)\lambda \Delta\lambda-m^2\lambda^2|\eta|^2\Big]\\
&&+2m^2\lambda^2|\eta|^2
+2(m-2)\left[|\nabla\lambda|^2+\lambda \Delta\lambda\right]\Big\}.
\end{eqnarray}
The proof of Proposition \ref{prop1} follows from equation \eqref{prop1-eq8}.
\end{proof}

\begin{corollary}
Let $\phi:(M^m,\bar{g})\to (N^n, h)$ $(m\neq4)$ be a conformal immersion with conformal factor $\lambda: M^m\to \r^+$.
If $m=2$ or $(m-2)|\nabla\lambda|^2=2\lambda \Delta\lambda$, then the condition $\bar S_2 = 0$ implies that the associated isometric immersion $\phi:(M^m,g)\to (N^n, h)$ is harmonic.
\end{corollary}

\begin{proof}
Observe that when $m=2$ or $(m-2)|\nabla\lambda|^2=2\lambda \Delta\lambda$,
Proposition \ref{prop1} yields
\begin{eqnarray}\label{co1-eq1}
{\rm Tr}_{\bar{g}}\,\bar{S}_2&=& -\frac{1}{2} m^2 (m-4)\lambda^4 |\eta|^2.
\end{eqnarray}
It follows from \eqref{co1-eq1} that $\bar S_2=0$ implies $\eta=0$. Hence $\tau(\phi,g)=0$.
\end{proof}

\begin{theorem}\label{th2}
Let $\phi:(M^m,\bar{g})\to (N^n, h)$ be a conformal immersion with conformal factor $\lambda: M^m\to \r^+$.
Then, the divergence of biharmonic stress tensor of the conformal immersion $\phi$ is given by
\begin{eqnarray*}
\lambda^{-2}\left[{\rm div}_{\bar{g}}\,\bar{S}_2\right](X)
&=& (m-2)\Big(\Delta\lambda-(m-2)\lambda^{-1}|\nabla\lambda|^2\Big)X(\lambda)+(m-2)\Big(\lambda X(\Delta\lambda)\\
& &-\frac{1}{2}(m-4)X(|\nabla\lambda|^2)+2\lambda {\rm Ric}^M(\nabla\lambda,X)\Big)
-\frac{1}{2}m^2\lambda^2 X(|\eta|^2)\\
& &+2m\lambda^2g({\rm Tr}_{g}\nabla^M A_\eta,X)
-2m(m-4)\lambda g(A_\eta(\nabla\lambda),X),
\end{eqnarray*}
where $A_\eta$ is the shape operator with respect to $\eta$ of the associated isometric immersion $\phi:(M^m,g)\to (N^n, h)$ with $\bar{g}=\lambda^{-2}g$.
\end{theorem}

\begin{proof}
Let $\{\bar{e}_i\}_{1\leq i\leq m}$ be a local orthonormal frame in $(M^m,\bar{g})$. We get $\{e_i=\frac{1}{\lambda}\bar{e}_i\}_{1\leq i\leq m}$ is a local orthonormal frame in $(M^m,g)$ with $g=\lambda^2\bar{g}$.
For all $X\in\Gamma(TM)$, we have
\begin{eqnarray}\label{th2-eq7}
\left[{\rm div}_{\bar{g}}\,\bar{S}_2\right](X)
&=&\nonumber\sum_{i=1}^m\left(\overline{\nabla}^M_{\bar e_i}\bar{S}_2\right)(X,\bar e_i)\\
&=&\nonumber\lambda^2\sum_{i=1}^m\left(\overline{\nabla}^M_{ e_i}\bar{S}_2\right)(X, e_i)\\
&=&\lambda^2\sum_{i=1}^m\left[ {e}_i\left(\bar{S}_2(X,{e}_i)\right)-\bar{S}_2(\overline{\nabla}_{{e}_i}^MX,{e}_i)
-\bar{S}_2(X,\overline{\nabla}_{{e}_i}^M{e}_i)\right],
\end{eqnarray}
where $\overline{\nabla}^M$ is the Levi-Civita of $(M^m,\bar{g})$. Using the identity $g(A_\eta(X), Y) = h(B(X, Y), \eta)$, we compute the first term of (\ref{th2-eq7}). For each $i = 1, \ldots, m$, we obtain
\begin{eqnarray}\label{th2-eq8}
\bar{S}_2(X,{e}_i)
&=&\nonumber \frac{1}{2}\Big[(m-2)(m-4)|\nabla\lambda|^2
-2(m-2)\lambda \Delta\lambda-m^2\lambda^2|\eta|^2\Big]g(X,e_i) \\
&& +2m\lambda^2g(A_{\eta}(X),e_i)+2(m-2)\Big[X(\lambda)e_i(\lambda)
+\lambda {\rm Hess\,}\lambda(X,e_i)\Big].
\end{eqnarray}
From the definition of ${\rm Hess\,}\lambda$, the symmetry of $A_{\eta}$, and equation (\ref{th2-eq8}), we get
\begin{eqnarray}\label{th2-eq9}
{e}_i\left(\bar{S}_2(X,{e}_i)\right)
&=&\nonumber e_i\Big( \frac{1}{2}\Big[(m-2)(m-4)|\nabla\lambda|^2
-2(m-2)\lambda \Delta\lambda-m^2\lambda^2|\eta|^2\Big]g(X,e_i)\\
&&\nonumber+2m\lambda^2g(A_{\eta}(e_i),X)+2(m-2)\left[g(X,\nabla\lambda)e_i(\lambda)+\lambda g(\nabla^M_{e_i}\nabla\lambda,X)\right]\Big)\\
&=&\nonumber\frac{1}{2}\Big[
(m-2)(m-4)e_i(|\nabla\lambda|^2)
-2(m-2)e_i(\lambda)\Delta\lambda
-2(m-2)\lambda e_i(\Delta\lambda)\\
&&\nonumber-2m^2\lambda e_i(\lambda)|\eta|^2
-m^2\lambda^2e_i(|\eta|^2)\Big]g(X,e_i)+\frac{1}{2}\Big[(m-2)(m-4)|\nabla\lambda|^2\\
&&\nonumber -2(m-2)\lambda \Delta\lambda-m^2\lambda^2|\eta|^2\Big]e_i(g(X,e_i))+4m\lambda e_i(\lambda)g(A_{\eta}(e_i),X)\\
&&\nonumber+2m\lambda^2g(\nabla_{e_i}^MA_{\eta}(e_i),X)+2m\lambda^2g(A_{\eta}(e_i),\nabla_{e_i}^MX)\\
&&\nonumber +2(m-2)g(\nabla_{e_i}^MX,\nabla\lambda)e_i(\lambda)+2(m-2)g(X,\nabla_{e_i}^M\nabla\lambda)e_i(\lambda) \\
&&\nonumber+2(m-2)g(X,\nabla\lambda)e_i(e_i(\lambda))+2(m-2)e_i(\lambda )g(\nabla^M_{e_i}\nabla\lambda,X)\\
&&+2(m-2)\lambda g(\nabla_{e_i}^M\nabla^M_{e_i}\nabla\lambda,X)+2(m-2)\lambda g(\nabla^M_{e_i}\nabla\lambda,\nabla_{e_i}^MX).
\end{eqnarray}
To compute the second term of (\ref{th2-eq7}), for all $i=1,...,m$, we have
\begin{eqnarray}\label{th2-eq10}
\bar{S}_2(\overline{\nabla}^M_{e_i}X,{e}_i)
&=&\nonumber \frac{1}{2}\Big[(m-2)(m-4)|\nabla\lambda|^2
-2(m-2)\lambda\Delta\lambda-m^2\lambda^2|\eta|^2\Big]g(\overline{\nabla}^M_{e_i}X,e_i)\\
&&\nonumber+2m\lambda^2g(A_{\eta}(e_i),\overline{\nabla}^M_{e_i}X)
+2(m-2)g(\overline{\nabla}^M_{e_i}X,\nabla\lambda)e_i(\lambda)\\
&&+2(m-2)\lambda g({\nabla}^M_{e_i}\nabla\lambda,\overline{\nabla}^M_{e_i}X).
\end{eqnarray}
By using the following formula
\begin{equation}\label{th2-eq11}
    \overline{\nabla}^M_{e_i}X={\nabla}^M_{e_i}X-\lambda^{-1}e_i(\lambda)X-\lambda^{-1}X(\lambda)e_i+\lambda^{-1}g(e_i,X)\nabla\lambda,
\end{equation}
the equation (\ref{th2-eq10}) becomes
\begin{eqnarray}\label{th2-eq12}
\bar{S}_2(\overline{\nabla}^M_{e_i}X,{e}_i)
&=&\nonumber \frac{1}{2}\Big[(m-2)(m-4)|\nabla\lambda|^2
-2(m-2)\lambda \Delta\lambda-m^2\lambda^2|\eta|^2\Big]\Big(g({\nabla}^M_{e_i}X,e_i)\\
&&\nonumber-\lambda^{-1}e_i(\lambda)g(X,e_i)-\lambda^{-1}X(\lambda)g(e_i,e_i)+\lambda^{-1}g(e_i,X)g(\nabla\lambda,e_i)\Big)\\
&&\nonumber+2m\lambda^2g(A_{\eta}(e_i),{\nabla}^M_{e_i}X)
-2m\lambda e_i(\lambda) g(A_{\eta}(e_i),X)
-2m\lambda X(\lambda) g(A_{\eta}(e_i),e_i)\\
&&\nonumber+2m\lambda g(e_i,X) g(A_{\eta}(e_i),\nabla\lambda)+2(m-2)g({\nabla}^M_{e_i}X,\nabla\lambda)e_i(\lambda)\\
&&\nonumber-2(m-2)\lambda^{-1}e_i(\lambda)g(X,\nabla\lambda)e_i(\lambda)-2(m-2)\lambda^{-1}X(\lambda)g(e_i,\nabla\lambda)e_i(\lambda)\\
&&\nonumber+2(m-2)\lambda^{-1}g(e_i,X)g(\nabla\lambda,\nabla\lambda)e_i(\lambda)
+2(m-2)\lambda g({\nabla}^M_{e_i}\nabla\lambda,{\nabla}^M_{e_i}X)\\
&&\nonumber -2(m-2) e_i(\lambda)g({\nabla}^M_{e_i}\nabla\lambda,X)-2(m-2) X(\lambda)g({\nabla}^M_{e_i}\nabla\lambda,e_i)\\
&&+2(m-2) g(e_i,X)g({\nabla}^M_{e_i}\nabla\lambda,\nabla\lambda).
\end{eqnarray}
To compute the last term of (\ref{th2-eq7}), for all $i=1,...,m$, we have
\begin{eqnarray}\label{th2-eq13}
\bar{S}_2(X,\overline{\nabla}^M_{e_i}{e}_i)
&=&\nonumber \frac{1}{2}\Big[(m-2)(m-4)|\nabla\lambda|^2
-2(m-2)\lambda\Delta\lambda-m^2\lambda^2|\eta|^2\Big]g(X,\overline{\nabla}^M_{e_i}e_i)\\
&&\nonumber+2m\lambda^2g(A_{\eta}(\overline{\nabla}^M_{e_i}e_i),X)
+2(m-2)g(X,\nabla\lambda)(\overline{\nabla}^M_{e_i}e_i)(\lambda)\\
&& +2(m-2)\lambda g({\nabla}^M_{\overline{\nabla}^M_{e_i}e_i}\nabla\lambda,X).
\end{eqnarray}
Note that for all $i=1,...,m$, we have
\begin{equation}\label{th2-eq14}
    \overline{\nabla}^M_{e_i}e_i={\nabla}^M_{e_i}e_i-2\lambda^{-1}e_i(\lambda)e_i+\lambda^{-1}g(e_i,e_i)\nabla\lambda.
\end{equation}
By using equations (\ref{th2-eq13}) and (\ref{th2-eq14}), we obtain
\begin{eqnarray}\label{th2-eq15}
\bar{S}_2(X,\overline{\nabla}^M_{e_i}{e}_i)
&=&\nonumber \frac{1}{2}\Big[(m-2)(m-4)|\nabla\lambda|^2
-2(m-2)\lambda \Delta\lambda-m^2\lambda^2|\eta|^2\Big]\Big(g(X,{\nabla}^M_{e_i}e_i)\\
&&\nonumber-2\lambda^{-1}e_i(\lambda)g(X,e_i)+\lambda^{-1}g(e_i,e_i)g(X,\nabla\lambda)\Big)
+2m\lambda^2g(A_{\eta}({\nabla}^M_{e_i}e_i),X)\\
&&\nonumber
    -4m\lambda e_i(\lambda)g(A_{\eta}(e_i),X)
    +2m\lambda g(e_i,e_i)g(A_{\eta}(\nabla\lambda),X)\\
    &&\nonumber+2(m-2)g(X,\nabla\lambda)({\nabla}^M_{e_i}e_i)(\lambda)
    -4(m-2)\lambda^{-1}e_i(\lambda)g(X,\nabla\lambda)e_i(\lambda)\\
    &&\nonumber+2(m-2)\lambda^{-1}g(e_i,e_i)g(X,\nabla\lambda)|\nabla\lambda|^2
    +2(m-2)\lambda g({\nabla}^M_{{\nabla}^M_{e_i}e_i}\nabla\lambda,X)\\
    &&-4(m-2) e_i(\lambda) g({\nabla}^M_{e_i}\nabla\lambda,X)
    +2(m-2)g(e_i,e_i) g({\nabla}^M_{\nabla\lambda}\nabla\lambda,X).
\end{eqnarray}
Theorem \ref{th2} follows from equations (\ref{th2-eq7}), (\ref{th2-eq9}), (\ref{th2-eq12}), (\ref{th2-eq15}), and the following
\begin{eqnarray*}
e_i(g(X,e_i))&=&g(\nabla^M_{e_i}X,e_i)+g(X,\nabla^M_{e_i}e_i),\quad\forall i=1,...,m,\\
\Delta \lambda &=& \sum_{i=1}^m\left[e_i(e_i(\lambda))-(\nabla^M_{e_i}e_i)(\lambda)\right],\\
{\rm Tr}_g({\nabla}^M)^2\nabla \lambda  &=& \sum_{i=1}^m\left[\nabla^M_{e_i}\nabla^M_{e_i}\nabla\lambda-\nabla^M_{\nabla^M_{e_i}e_i}\nabla\lambda\right]\\
  &=& {\rm Ric}^M\,(\nabla\lambda)+\nabla {\Delta}\lambda,\\
{\nabla}^M_{\nabla\lambda}\nabla\lambda &=&\frac{1}{2} \nabla(|\nabla\lambda|^2),\\
|\eta|^2&=&\frac{1}{m}\sum_{i=1}^mg(A_\eta(e_i),e_i).
\end{eqnarray*}
\end{proof}

Recall that a \textbf{biconservative} map is a smooth map between Riemannian manifolds whose stress-energy tensor associated with the bienergy functional is divergence-free. For a conformal immersion, we obtain the following theorem.

\begin{theorem}
Let $\phi:(M^m,\bar{g})\to (N^n, h)$ be a conformal immersion with conformal factor $\lambda: M^m\to \r^+$.
Then, $\phi$ is biconservative conformal immersion if and only if
\begin{eqnarray*}
& & (m-2)\Big(\Delta\lambda-(m-2)\lambda^{-1}|\nabla\lambda|^2\Big)X(\lambda)+(m-2)\Big(\lambda X(\Delta\lambda)\\
& &-\frac{1}{2}(m-4)X(|\nabla\lambda|^2)+2\lambda {\rm Ric}^M(\nabla\lambda,X)\Big)
-\frac{1}{2}m^2\lambda^2 X(|\eta|^2)\\
& &+2m\lambda^2g({\rm Tr}_{g}\nabla^M A_\eta,X)
-2m(m-4)\lambda g(A_\eta(\nabla\lambda),X)=0.
\end{eqnarray*}
\end{theorem}

%%%%%%%%%%%%%%%%%%%%%%%%%%%%%%%%%%%%%%%%%%%%%%%%%%%%%%%%%%  Hypersurfaces Case %%%%%%%%%%%%%%%%%%%%%%%%%%%%%%%%%%%%%%%%%%%%%%%%%%%%

\section{The stress-bienergy tensor of conformal hypersurfaces}

We specialize to conformal hypersurfaces $\phi:(M^m,\bar g)\to(N^{m+1},h)$. Let
$\xi$ be a unit normal vector field of the associated isometric hypersurface
$\phi:(M^m,g)\to(N^{m+1},h)$, with mean curvature function $H$ and shape operator
$A$, so that $\eta=H\xi$ and $A_\eta=HA$. The formulas of the previous section then
simplify, and we show that biconservativity is equivalent to the vanishing of the
tangent part of the bitension field.

\begin{corollary}\label{co2}
Let $\phi:(M^m,\bar{g})\to (N^{m+1}, h)$ be a conformal hypersurface with conformal factor $\lambda: M^m\to \r^+$.
Then, the biharmonic stress tensor of the conformal hypersurface $\phi$ is given by
\begin{eqnarray*}
 \bar{S}_2(X,Y)
   &=& \frac{1}{2}\Big[(m-2)(m-4)|\nabla\lambda|^2
       -2(m-2)\lambda \Delta\lambda-m^2H^2\lambda^2\Big]g(X,Y)\\
   &&  +2mH\lambda^2g(A(X),Y)
       +2(m-2)\left[X(\lambda)Y(\lambda)+\lambda ({\rm Hess\,}\lambda)(X,Y)\right],
\end{eqnarray*}
for all $X,Y\in\Gamma(TM)$, where $H$ and $A$ are the mean curvature function and the shape operator with respect to unit normal vector field $\xi$ of the associated hypersurface $\phi:(M^m,g)\to (N^{m+1}, h)$, respectively.
Moreover, if $(M^m,g)$ is a totally umbilical hypersurface in $(N^{m+1}, h)$, then the biharmonic stress tensor of the conformal hypersurface $\phi$ is given by
\begin{eqnarray*}
 \bar{S}_2(X,Y)
   &=& \frac{1}{2}\Big[(m-2)(m-4)|\nabla\lambda|^2
       -2(m-2)\lambda \Delta\lambda
       -m(m-4)H^2\lambda^2
       \Big]g(X,Y)\\
   & & +2(m-2)\left[X(\lambda)Y(\lambda)+\lambda ({\rm Hess\,}\lambda)(X,Y)\right].
\end{eqnarray*}
\end{corollary}

\begin{proof}
Using Theorem \ref{th1} with $\eta=H\xi$, and hence $A_\eta=H A$, together with the identity $g(A_\eta(X), Y) = h(B(X,Y), \eta)$.
If $(M^m,g)$ is a totally umbilical hypersurface in $(N^{m+1}, h)$ we have the shape operator $A = H Id$, where $Id$ is the
identity map.
\end{proof}

\begin{corollary}\label{co3}
Let $\phi:(M^m,\bar{g})\to (N^{m+1}, h)$ be a conformal hypersurface with conformal factor $\lambda: M^m\to \r^+$.
Then, the divergence of biharmonic stress tensor of the conformal hypersurface $\phi$ is given by
\begin{eqnarray*}
\lambda^{-2}\left[{\rm div}_{\bar{g}}\,\bar{S}_2\right](X)
   &=& (m-2)\Big(\Delta\lambda-(m-2)\lambda^{-1}|\nabla\lambda|^2\Big)X(\lambda)+(m-2)\Big(\lambda X(\Delta\lambda)\\
   & &-\frac{1}{2}(m-4)X(|\nabla\lambda|^2)+2\lambda {\rm Ric}^M(\nabla\lambda,X)\Big)
     -2mH\lambda^2{\rm Ric}^N(\xi,X)\\
   & &+m^2H\lambda^2X(H)+2m\lambda^2g(A(\nabla H),X)-2m(m-4)H\lambda g(A(\nabla\lambda),X).
\end{eqnarray*}
Moreover, if $(M^m,g)$ is a totally umbilical hypersurface in a space form $(N^{m+1}(c), h)$ $(m\geq3)$, then the biharmonic stress tensor of the conformal hypersurface $\phi$ is given by
\begin{eqnarray*}
\lambda^{-2}\left[{\rm div}_{\bar{g}}\,\bar{S}_2\right](X)
   &=& \Big((m-2)\Delta\lambda-(m-2)^2\lambda^{-1}|\nabla\lambda|^2-2m(m-4)H^2\lambda\\
   & & +2(m-1)(m-2)(c+H^2)\lambda\Big)X(\lambda)+(m-2)\lambda X(\Delta\lambda)\\
   & & -\frac{1}{2}(m-2)(m-4)X(|\nabla\lambda|^2).
\end{eqnarray*}
\end{corollary}

\begin{proof}
Using Theorem \ref{th2} together with $\eta=H\xi$, $A_\eta=H A$, and the identity
${\rm Tr}_{g}\nabla^M A=m\,\nabla H-({\rm Ric}^N\xi)^\top$ (see \cite{MM}). When
$(M^m,g)$ is a totally umbilical hypersurface in $(N^{m+1}(c), h)$ using $A=H Id$. By using Gauss's equation we find that
${\rm Ric}^M\,(X)=(m-1)(H^2+c)X$ for all $X\in\Gamma(TM)$. Hence, for $m\geq3$ the curvature function $H$ is a constant.
\end{proof}

\begin{theorem}\label{th3}
Let $\phi:(M^m,\bar{g})\to (N^{m+1}, h)$ be a conformal hypersurface with conformal factor $\lambda: M^m\to \r^+$.
Then, ${\rm div}_{\bar{g}}\,\bar{S}_2=0$ if and only if the tangent part of $\tau_{2}(\phi, \bar g)$ vanishes.
\end{theorem}

\begin{proof}
We know that the tangential component of the bitension field $\tau_{2}(\phi, \bar g)$ vanishes if and only if
\begin{align}\label{th3-eq20}
&-m\Big[2A(\nabla H)+\frac{m}{2}\nabla H^2-2H({\rm Ric}^N\xi)^\top\Big]-(m-2)\Big[2\, {\rm Ric}^M\,(\nabla \ln\lambda)\\\notag
&+\nabla {\Delta}\ln\lambda\Big]-2(m-2)\Big[\Delta {\rm ln}\lambda-(m-4)\left|\nabla\ln\lambda\right|^2\Big]\nabla\ln\lambda\\\notag
&+\frac{1}{2}(m-2)(m-6)\nabla |\nabla\ln\lambda|^2+2m(m-4)HA(\nabla\ln\lambda)=0,
\end{align}
(see \cite{MO}). By using $\Delta \ln \lambda=\lambda^{-1}\Delta \lambda-\lambda^{-2}|\nabla\lambda|^2, \, |\nabla\ln\lambda|^2=\lambda^{-2}|\nabla\lambda|^2$ and a straightforward calculation, equation (\ref{th3-eq20}) becomes
\begin{align*}
&-m\Big[2A(\nabla H)+\frac{m}{2}\nabla H^2-2H({\rm Ric}^N\xi)^\top\Big]
 -(m-2)\Big[2\lambda^{-1} {\rm Ric}^M\,(\nabla\lambda)-\lambda^{-2}({\Delta}\lambda)\nabla\lambda\\
&+\lambda^{-1}\nabla\Delta\lambda+2\lambda^{-3}|\nabla\lambda|^2\nabla\lambda
 -\lambda^{-2}\nabla|\nabla\lambda|^2\Big]
 -2(m-2)\Big[\lambda^{-1}\Delta \lambda-\lambda^{-2}|\nabla\lambda|^2\\
&-(m-4)\lambda^{-2}|\nabla\lambda|^2\Big]\lambda^{-1}\nabla\lambda
 +\frac{1}{2}(m-2)(m-6)\Big[-2\lambda^{-3}|\nabla\lambda|^2\nabla\lambda+\lambda^{-2}\nabla|\nabla\lambda|^2\Big]\\
&+2m(m-4)H\lambda^{-1}A(\nabla\lambda)=0,
\end{align*}
this is equivalent to the following
\begin{align}\label{th3-eq21}
&-(m-2)\lambda^{-2}(\Delta\lambda)\nabla\lambda
 +(m-2)^2\lambda^{-3}|\nabla\lambda|^2\nabla\lambda
 -(m-2)\lambda^{-1} \nabla\Delta\lambda\\\notag
&+\frac{1}{2}(m-2)(m-4)\lambda^{-2}\nabla|\nabla\lambda|^2
 -2(m-2)\lambda^{-1}{\rm Ric}^M\,(\nabla\lambda)
 +2mH({\rm Ric}^N\xi)^\top\\\notag
&-m^2H\nabla H
 -2mA(\nabla H)
 +2m(m-4)H\lambda^{-1}A(\nabla\lambda)=0.
\end{align}
By multiplying (\ref{th3-eq21}) by $\lambda^2$ and applying it to $X \in \Gamma(TM)$, we obtain that
\begin{align}\label{th3-eq22}
&-(m-2)(\Delta\lambda)X(\lambda)
 +(m-2)^2\lambda^{-1}|\nabla\lambda|^2X(\lambda)
 -(m-2)\lambda X(\Delta\lambda)\\\notag
&+\frac{1}{2}(m-2)(m-4)X(|\nabla\lambda|^2)
 -2(m-2)\lambda g({\rm Ric}^M\,(\nabla\lambda),X)
 +2mH\lambda^2h({\rm Ric}^N\xi,X)\\\notag
&-m^2\lambda^2HX(H)
 -2m\lambda^2g(A(\nabla H),X)
 +2m(m-4)H\lambda g(A(\nabla\lambda),X)=0.
\end{align}
Note that (\ref{th3-eq22}) is equivalent to $[{\rm div}_{\bar{g}}\,\bar{S}_2](X)=0$.
\end{proof}

\begin{remark}
The equivalence in Theorem \ref{th3} is in fact a general feature of immersions. By \cite{LMO}, for any smooth map $\psi:(M,g)\to(N,h)$,
$\operatorname{div} S_2=-\langle \tau_2(\psi),d\psi\rangle$. If $\psi$ is an immersion, the vectors $d\psi(X)$, $X\in\Gamma(TM)$, span the tangent space of $\psi(M)$. Hence $\operatorname{div}S_2=0$ if and only if the tangent part of $\tau_2(\psi)$ vanishes. Applying this to the conformal immersion $\phi:(M^m,\bar g)\to(N^n,h)$, we obtain, in any codimension, $\operatorname{div}_{\bar g}\bar S_2=0$ if and only if $\big(\tau_2(\phi,\bar g)\big)^{\top}=0$.
Thus a conformal immersion is biconservative if and only if the tangent part of its bitension field vanishes. For conformal hypersurfaces, Theorem \ref{th3} recovers this fact by a direct computation, based on the expression \eqref{th3-eq20} of $\big(\tau_2(\phi,\bar g)\big)^{\top}$ from \cite{MO}.
\end{remark}

%%%%%%%%%%%%%%%%%%%%%%%%%%%%%%%%%%%%%%%%%%%%%%%%%%%%%%%%%   30/05/2026  Cherif %%%%%%%%%%%%%%%%%%%%%%%%%%%%%%%%%%%%%%%%%%%%%%%%%%%%%%%%%%%%%%%%%%

\section{Conformally biconservative hypersurfaces}

Given a hypersurface $(M^m,g)$ in $(N^{m+1},h)$, we ask for which conformal
factors $\lambda$ the immersion $\phi:(M^m,\lambda^{-2}g)\to(N^{m+1},h)$ is
biconservative. For totally umbilical hypersurfaces of space forms, we show that a
non-constant $\lambda$ gives a biconservative immersion if and only if it is an
isoparametric function whose structure functions satisfy an ODE, and we then give
explicit examples.

\begin{definition}
A hypersurface $(M^m,g)$ in a Riemannian manifold $(N^{m+1},h)$ is said to be conformally biconservative if there exists a smooth function $\lambda: M^m\to \r^+$ such that the conformal immersion $\phi : (M^m,\bar{g}=\lambda^{-2}g) \to (N^{m+1}, h)$ is biconservative, i.e., $\operatorname{div}_{\bar{g}} \bar{S}_2 = 0$.
\end{definition}

According to Corollary \ref{co3}, we conclude that

\begin{proposition}\label{prop2}
A totally umbilical hypersurface $(M^m,g)$ $(m\geq3)$ in a space form $(N^{m+1}(c),h)$ is conformally biconservative with conformal factor $\lambda: M^m\to \r^+$
if and only if
\begin{eqnarray*}
   & & \Big((m-2)\Delta\lambda-(m-2)^2\lambda^{-1}|\nabla\lambda|^2-2m(m-4)H^2\lambda\\
   & & +2(m-1)(m-2)(c+H^2)\lambda\Big)\nabla\lambda+(m-2)\lambda \nabla\Delta\lambda\\
   & & -\frac{1}{2}(m-2)(m-4)\nabla|\nabla\lambda|^2=0.
\end{eqnarray*}
\end{proposition}

Recall that a smooth function $f : M \to \mathbb{R}$ is said to be \textbf{isoparametric} if there exist smooth functions $\alpha,\beta : \mathbb{R} \to \mathbb{R}$ satisfying $|\nabla f|^2 = \alpha(f)$ and $\Delta f = \beta(f)$. From Proposition \ref{prop2}, we get the following Theorem.

\begin{theorem}\label{thm3}
A totally umbilical hypersurface $(M^m,g)$ $(m\geq3)$ in a space form $(N^{m+1}(c),h)$ is conformally biconservative with a non-constant conformal factor $\lambda: M^m\to \r^+$ if and only if $\lambda$ is an isoparametric function on $(M^m,g)$ with $|\nabla \lambda|^2=\alpha(\lambda)$ and $\Delta\lambda=\beta(\lambda)$ for some smooth real functions $\alpha$ and $\beta$ solving the ODE
\begin{eqnarray}\label{eq23}
   & &\nonumber -\frac{1}{2}(m-4)\alpha'(\lambda)-(m-2)\lambda^{-1}\alpha(\lambda)+\lambda \beta'(\lambda)+\beta(\lambda)-\frac{2m(m-4)}{m-2}H^2\lambda\\
   & &+2(m-1)(c+H^2)\lambda =0.
\end{eqnarray}
\end{theorem}

\begin{proof}
It follows from Proposition \ref{prop2} that a totally umbilical hypersurface $(M^m,g)$ in a space form $(N^{m+1}(c),h)$ can be conformally biconservative
if and only if
\begin{eqnarray}\label{eq24}
\nabla\Big[\lambda\Delta\lambda-\frac{1}{2}(m-4)|\nabla\lambda|^2\Big]
    &=&\nonumber\Big[(m-2)\lambda^{-1}|\nabla\lambda|^2+\frac{2m(m-4)}{m-2}H^2\lambda\\
    & &-2(m-1)(c+H^2)\lambda\Big]\nabla\lambda.
\end{eqnarray}
By using this and Lemma 2 of \cite{BO} we have
\begin{equation}\label{eq25}
    \lambda\Delta\lambda-\frac{1}{2}(m-4)|\nabla\lambda|^2=u(\lambda),
\end{equation}
for some smooth function $u$ with
\begin{equation}\label{eq26}
    u'(\lambda)=(m-2)\lambda^{-1}|\nabla\lambda|^2+\frac{2m(m-4)}{m-2}H^2\lambda-2(m-1)(c+H^2)\lambda.
\end{equation}
It follows from (\ref{eq25}) and (\ref{eq26}) that $\lambda$ is an isoparametric function. By substituting
$|\nabla \lambda|^2=\alpha(\lambda)$ and $\Delta\lambda=\beta(\lambda)$ into (\ref{eq24}) we get (\ref{eq23}).
\end{proof}

\begin{example}
A part $M^m$ of the sphere $\mathbb{S}^m$ can be conformally biconservative in $\mathbb{R}^{m+1}$ $(m\geq 3)$ with conformal factor
$\lambda=(1+|x|^2)^{-1}$, where $x=(x_1,\ldots,x_m)$ denotes the conformal coordinates on $\mathbb{S}^m$, if and only if $m=4$.\\
Indeed; by using the stereographic projection the part $M^m$ can be expressed as
\begin{eqnarray*}
 \phi:(\r^m,g_{ij}=\frac{4\delta_{ij}}{(1+|x|^2)^2})  &\rightarrow& (\mathbb{R}^{m+1},h_{uv}=\delta_{uv}). \\
   (x_1,...,x_m)&\mapsto& \left(\frac{2x_1}{1+|x|^2},...,\frac{2x_m}{1+|x|^2},\frac{-1+|x|^2}{1+|x|^2}\right)
\end{eqnarray*}
Now, we apply Theorem \ref{thm3} with $\lambda(x)=(1+|x|^2)^{-1}$, $H=1$, and $c=0$. A direct calculation gives
\begin{eqnarray*}
  \nabla\lambda &=& -\frac{1}{2}\,x_i\frac{\partial}{\partial x_i}, \\[0.3em]
  |\nabla\lambda|^2 &=& \frac{|x|^2}{(1+|x|^2)^{2}}, \\[0.3em]
  \Delta\lambda &=& \frac{m}{2}\frac{-1+|x|^2}{1+|x|^2}.
\end{eqnarray*}
A further computation yields
\begin{eqnarray*}
\alpha(\lambda)&=&|\nabla\lambda|^2 = \lambda-\lambda^2, \\[0.3em]
\beta(\lambda)&=&\Delta\lambda = \frac{m}{2}(1-2\lambda).
\end{eqnarray*}
Substituting these expressions into (\ref{eq23}), we obtain
$$(m-4)\left(m-2+4\lambda\right)=0.$$
Hence, $m=4$. Therefore, any open part of $\mathbb{S}^4$ is conformally biconservative in $\mathbb{R}^{5}$.
\end{example}

\begin{example}
Let $a,b \in \mathbb{R}$ with $a>0$ and $a^2+b^2=1$. For $m\geq 3$, consider the totally umbilical hypersurface
$\phi : \mathbb{S}^m(a) \to \mathbb{S}^{m+1},$ $\phi(v)=(v,b).$
A part $M^m$ of $\mathbb{S}^m(a)$ is conformally biconservative with conformal factor $\lambda=(1+|x|^2)^{-1},$
where $x=(x_1,...,x_m)$ denotes the conformal coordinates on $\mathbb{S}^m(a)$, if and only if $m=4$. \\
Indeed, in conformal coordinates induced by stereographic projection, $M^m$ can be written as
\begin{eqnarray*}
 \phi:(\r^m,g=\frac{4a^2\delta_{ij}}{(1+|x|^2)^2})  &\rightarrow& (\mathbb{S}^{m+1},h_{uv}=\frac{4\delta_{uv}}{(1+|y|^2)^2}). \\
   (x_1,...,x_m)&\mapsto& \frac{a}{1-b}\left(\frac{2x_1}{1+|x|^2},...,\frac{2x_m}{1+|x|^2},\frac{-1+|x|^2}{1+|x|^2}\right)
\end{eqnarray*}
We apply Theorem \ref{thm3} with $\lambda(x)=(1+|x|^2)^{-1}$, $H^2=b^2/(1-b^2)$, and $c=1$. A direct calculation gives
\begin{eqnarray*}
  \nabla\lambda &=& -\frac{1}{2a^2}\,x_i\frac{\partial}{\partial x_i}, \\[0.3em]
  |\nabla\lambda|^2 &=& \frac{1}{a^2}\,\frac{|x|^2}{(1+|x|^2)^{2}}, \\[0.3em]
  \Delta\lambda   &=&  \frac{m}{2a^2}\frac{-1+|x|^2}{1+|x|^2}.
\end{eqnarray*}
A further computation yields
\begin{eqnarray*}
\alpha(\lambda) &=& |\nabla\lambda|^2 =\frac{1}{a^2}(\lambda-\lambda^2), \\[0.3em]
\beta(\lambda) &=& \Delta\lambda =  \frac{m}{2a^2}(1-2\lambda).
\end{eqnarray*}
Substituting these expressions into (\ref{eq23}), we obtain
$$(m-4)\left(m-2+2[m(b-1)(b+1)+2]\lambda\right)=0.$$
Hence, $m=4$ and any open part of $\mathbb{S}^4(a)$ is conformally biconservative in $\mathbb{S}^{5}$.
\end{example}

\begin{example}
A part $M^m$ $(m\geq 3)$ of the totally umbilical hypersurface $$\phi : (\mathbb{R}^m\backslash\{0\},g_{ij}=\delta_{ij}) \to (\mathbb{H}^{m+1},h_{uv}=\frac{\delta_{uv}}{y_{m+1}^2}),\quad \phi(x)=(x,1),$$ is conformally biconservative with conformal factor $\lambda=|x|^2$ if and only if $m=4$. \\
Indeed, we apply Theorem \ref{thm3} with $H=1$, and $c=-1$. We have
\begin{eqnarray*}
  \nabla\lambda &=& 2 x_i \frac{\partial}{\partial x_i}, \\[0.3em]
  |\nabla\lambda|^2 &=& 4 |x|^2, \\[0.3em]
  \Delta\lambda   &=& 2m.
\end{eqnarray*}
Therefore, the functions $\alpha$ and $\beta$ are given by
\begin{eqnarray*}
\alpha(\lambda) &=& |\nabla\lambda|^2 =4\lambda, \\[0.3em]
\beta(\lambda) &=& \Delta\lambda = 2m.
\end{eqnarray*}
Substituting these expressions into (\ref{eq23}), we obtain
$$(m-4)\left(m-2+\frac{m}{2}\lambda\right)=0.$$
Since $\lambda=|x|^2>0$ and $m\geq 3$, we have $m-2+\frac m2\lambda>0$, hence $m=4$.
\end{example}

\begin{example}
Let \(a_1,...,a_m,b,s\in\mathbb{R}\) with \(s\neq 0\), and assume that $a_ix_i+b>1$
on an open set \(\Omega\subset\mathbb{R}^m\) $(m\geq 3)$. Consider the totally umbilical hypersurface
\begin{eqnarray*}
  \phi: (\Omega,g) &\rightarrow& (\mathbb{R}^{m+1},h=dy_1^2+e^{2sy_1}(dy_2^2+...+dy_{m+1}^2)).  \\
   (x_1,...,x_m) &\mapsto& (-\frac{1}{s}\ln(a_i x_i+b),x_1,...,x_m)
\end{eqnarray*}
Then, \(\phi\) is a conformally biconservative hypersurface with
$\lambda=\bigl[\ln\bigl(a_ix_i+b\bigr)\bigr]^{1/2}$ if and only if \(m=3\) and $s=\pm\sqrt{\frac{2}{3}}\,|a|$ where $|a|^2=a_1^2+...+a_m^2>0$.\\
Indeed, the components of the induced Riemannian metric
$g = \phi^* h$ of this hypersurface are given by $g_{ij}=(s^{-2}a_ia_j+\delta_{ij})/(a_kx_k+b)^2$.
We apply Theorem \ref{thm3} with $H^2=s^4/(|a|^2+s^2)$ and $c=-s^2$. We have
\begin{eqnarray*}
  \nabla\lambda &=& \frac{s^2}{2}\frac{a_i(a_jx_j+b)}{(|a|^2+s^2)\bigl[\ln\bigl(a_kx_k+b\bigr)\bigr]^{1/2}} \frac{\partial}{\partial x_i}, \\[0.3em]
  |\nabla\lambda|^2 &=& \frac{s^2}{4} \frac{|a|^2}{(|a|^2+s^2)\ln\bigl(a_ix_i+b\bigr)}, \\[0.3em]
  \Delta\lambda   &=& -\frac{s^2}{4}\frac{|a|^2\bigl[1+2(m-1)\ln(a_ix_i+b)\bigr]}{(|a|^2+s^2)\bigl[\ln\bigl(a_jx_j+b\bigr)\bigr]^{3/2}}.
\end{eqnarray*}
So that, the functions $\alpha$ and $\beta$ are given by
\begin{eqnarray*}
\alpha(\lambda) &=& |\nabla\lambda|^2 = \frac{s^2}{4} \frac{|a|^2}{(|a|^2+s^2)\lambda^2}, \\[0.3em]
\beta(\lambda) &=& \Delta\lambda = -\frac{s^2}{4}\frac{|a|^2\bigl[1+2(m-1)\lambda^2\bigr]}{(|a|^2+s^2)\lambda^3}.
\end{eqnarray*}
Substituting the above expressions into \eqref{eq23}, we obtain
\[
m(m-4)s^2+(m-1)(m-2) |a|^2=0.
\]
Consequently, a real solution exists only when \(m=3\), in which case
$s=\pm\sqrt{\frac{2}{3}}\,|a|$.
\end{example}

%%%%%%%%%%%%%%%%%%%%%%%%%%%%%%%%%%%%%%%%%%%%%%%%%%%%%%%%%%%%%%%%%%   31/05/2026 Cherif  %%%%%%%%%%%%%%%%%%%%%%%%%%%%%%%%%%%%%%%%%%%%%%%%%%%%%%%%

\section{Conformal immersions with vanishing stress-bienergy tensor}

Finally, we study the equation $\bar S_2=0$. We first write it in a form that
involves only $\lambda$, $B$ and $\eta$. We then show that, on a compact orientable
manifold of dimension $m\neq 4$, it forces the associated isometric immersion to
be harmonic, and we give a non-compact example with non-constant $\lambda$.

\begin{proposition}\label{prop3}
Let $\phi:(M^m,\bar{g})\to (N^n, h)$ be a conformal immersion with conformal factor $\lambda: M^m\to \r^+$.
Assume that $m\neq2$. Then, $\bar{S}_2=0$ if and only if
\begin{eqnarray}\label{prop3-eq27}
   & &\nonumber\frac{1}{m-2}\bigl[m^2\lambda^2|\eta|^2+(m-2)^2|\nabla\lambda|^2\bigr]g(X,Y)-2m\lambda^2h(B(X,Y),\eta )\\
   & &-2(m-2)\bigl[X(\lambda) Y(\lambda)+\lambda ({\rm Hess\,}\lambda)(X,Y)\bigr] =0,
\end{eqnarray}
for all $X,Y\in\Gamma(TM)$.
\end{proposition}

\begin{proof}
Using Proposition 3.3 of \cite{LMO}, it follows that $\bar{S}_2=0$ if and only if
\begin{eqnarray}\label{prop3-eq28}
   & &\frac{1}{m-2}|\tau(\phi,\bar{g})|^2\bar{g}(X,Y)
    +h( {\rm d}\phi(X),\nabla_Y^\phi \tau(\phi,\bar{g}))+ h( {\rm d}\phi(Y),\nabla_X^\phi \tau(\phi,\bar{g}))=0,
\end{eqnarray}
for all $X,Y\in\Gamma(TM)$. By substituting (\ref{th1-eq3}) and (\ref{th1-eq4}) into (\ref{prop3-eq28}), we get (\ref{prop3-eq27}).
\end{proof}

\begin{corollary}\label{prop3}
Let $\phi:(M^m,\bar{g})\to (N^{m+1}, h)$ be a conformal hypersurface with conformal factor $\lambda: M^m\to \r^+$.
Assume that $m\neq2$. Then, $\bar{S}_2=0$ if and only if
\begin{eqnarray}\label{prop3-eq29}
   & &\nonumber\frac{1}{m-2}\bigl[m^2H^2\lambda^2+(m-2)^2|\nabla\lambda|^2\bigr]g(X,Y)-2mH\lambda^2g(A(X),Y)\\
   & &-2(m-2)\bigl[X(\lambda) Y(\lambda)+\lambda ({\rm Hess\,}\lambda)(X,Y)\bigr] =0,
\end{eqnarray}
for all $X,Y\in\Gamma(TM)$. Moreover, if $(M^m,g)$ is a totally umbilical hypersurface in $(N^{m+1},h)$,
then $\bar{S}_2=0$ if and only if
\begin{eqnarray}\label{prop3-eq30}
{\rm Hess\,}\lambda^2=2(\lambda\,{\rm Hess\,}\lambda+{\rm d}\lambda\otimes{\rm d}\lambda)=\Big[|\nabla\lambda|^2-\frac{m(m-4)}{(m-2)^2}H^2\lambda^2\Big]g.
\end{eqnarray}
\end{corollary}

\begin{theorem}\label{thm:compact}
Let $\phi:(M^m,\bar g)\to(N^n,h)$, with $m\neq 4$, be a conformal immersion
with conformal factor $\lambda:M^m\to\mathbb{R}^+$. Assume that $\bar S_2=0$
and that $M$ is compact and orientable. Then the associated isometric
immersion $\phi:(M^m,g)\to(N^n,h)$, $g=\lambda^2\bar g$, is harmonic.
If moreover $m\neq 2$, then $\lambda$ is constant.
\end{theorem}

\begin{proof}
Theorem 3.6 of \cite{LMO} states that a map $\psi:(M^m,\bar g)\to(N^n,h)$ from a
compact orientable Riemannian manifold with $m\neq 4$ and vanishing
stress-bienergy tensor is harmonic. Applying it to
$\phi:(M,\bar g)\to(N,h)$, the hypothesis $\bar S_2=0$ gives $\tau(\phi,\bar g)=0$.
By equation (\ref{th1-eq2}), $\tau(\phi,\bar g)=m\lambda^2\eta-(m-2)\lambda\, d\phi(\nabla\lambda)$,
where the first term is normal and the second is tangent to $\phi(M)$. Hence,
$\eta=0$ and $(m-2)\,d\phi(\nabla\lambda)=0 $. The first equality gives $\tau(\phi,g)=m\eta=0$, so
$\phi:(M,g)\to(N,h)$ is harmonic. If $m\neq 2$, the second gives
$d\phi(\nabla\lambda)=0$, hence $\nabla\lambda=0$ because $\phi$ is an
immersion. Thus, $\lambda$ is constant on each connected component of $M$.
\end{proof}

\begin{remark}
For $m=2$, Remark \ref{remark1}\,(ii) gives $\bar S_2=\lambda^2S_2$, where $S_2$
is the stress-bienergy tensor of the isometric immersion $\phi:(M^2,g)\to(N^n,h)$. By (7),
$S_2(X,Y)=-\tfrac12|\tau(\phi,g)|^2g(X,Y)+2h(\tau(\phi,g),B(X,Y))$,
so $S_2=0$ whenever $\phi$ is minimal, since then $\tau(\phi,g)=2\eta=0$.
Hence, every compact orientable minimal surface $\phi:(M^2,g)\to(N^n,h)$ (for
instance a minimal torus in $\mathbb{S}^3$), taken with an arbitrary smooth function
$\lambda:M^2\to\mathbb{R}^+$ and $\bar g=\lambda^{-2}g$, satisfies
$\bar S_2=0$ with $\lambda$ non-constant.
\end{remark}

\begin{example}
For $a,b,s,t\in\mathbb{R}$ and $m\geq 3$, the totally umbilical hypersurface
\[
\phi:\Big(\mathbb{H}^m,\ g_{ij}=\frac{\delta_{ij}\,(1+\delta_{im}\delta_{jm}a^2)}{x_m^{s}}\Big)\to
\Big(\mathbb{H}^{m+1},\ h_{uv}=\frac{\delta_{uv}}{y_{m+1}^{s}}\Big),\quad
\phi(x)=(ax_m+b,x),
\]
is conformally vanishing biharmonic stress tensor with conformal factor
$\lambda=x_m^{t}$, i.e., $\bar S_2=0$ if and only if
\begin{equation}\label{31}
-st=2t\Big(2t-1+\frac s2\Big)=t^2-\frac{m(m-4)a^2s^2}{4(m-2)^2}.
\end{equation}
Indeed, we apply equation (\ref{prop3-eq30}) with
$H^2=\dfrac{a^2s^2x_m^{s-2}}{4(1+a^2)}$. A direct calculation gives
\[
\nabla\lambda=\frac{t}{1+a^2}\,x_m^{s+t-1}\frac{\partial}{\partial x_m},\quad
|\nabla\lambda|^2=\frac{t^2}{1+a^2}\,x_m^{s+2t-2},
\]
\[
\operatorname{Hess}\lambda^2=x_m^{2(t-1)}
\begin{pmatrix}
-\dfrac{st}{1+a^2} & 0 & \cdots & 0\\
0 & -\dfrac{st}{1+a^2} & \cdots & 0\\
\vdots & \vdots & \ddots & \vdots\\
0 & 0 & \cdots & 2t\big(2t-1+\frac s2\big)
\end{pmatrix}.
\]
Substituting these into (\ref{prop3-eq30}), we deduce equations \eqref{31}. Observe that, for
$m=4$, we get $s=-t=-1$ or $t=0$. For $m\neq 4$ and $t\neq \tfrac12$, we get
$s=1-2t$ and
\[
a=\pm\frac{2(m-2)}{1-2t}\sqrt{\frac{t(1-t)}{m(m-4)}},
\]
which is real if and only if $t\in[0,1]\backslash\{\tfrac12\}$ when $m\geq 5$, and
$t\in(-\infty,0]\cup[1,\infty)$ when $m=3$.
\end{example}

\section{Concluding remarks}

We computed the stress-bienergy tensor of a conformal immersion, together with its
trace and its divergence, in terms of the conformal factor and of the second
fundamental form of the associated isometric immersion. For totally umbilical
hypersurfaces of space forms, conformal biconservativity with a non-constant
conformal factor is equivalent to $\lambda$ being isoparametric, with structure
functions $(\alpha,\beta)$ solving the ODE \eqref{eq23}. The examples show that this
occurs in dimension $m=4$ for spheres, small spheres and horospheres, and in
dimension $m=3$ for a hypersurface of a warped product. When $\bar S_2=0$, $M$ is
compact and $m\neq 4$, the associated isometric immersion is harmonic, whereas non-compact examples with non-constant $\lambda$
exist.

\subsection*{Acknowledgments} The author would like to thank Professor Y.-L. Ou for his useful remarks and suggestions.

\subsection*{Funding Declaration} No funding was received to assist with the preparation of this manuscript.

\subsection*{Competing Interest declaration} The author has no competing interests to declare that are relevant to the content of this article.

\subsection*{Data availability} Not applicable.

\end{document}